# Quasi-linear analysis of dispersion relation preservation for nonlinear schemes

Fengyuan Xu [1], Pan Yan [1], Qin Li [1*], Yancheng You [1]

[1]School of Aerospace Engineering, Xiang'an district, Xiamen University, Xiamen, 361102, China

Email address:

    Fengyuan Xu: 35020191151118@stu.xmu.edu.cn

    Pan Yan: 35020201151395@stu.xmu.edu.cn

    Qin Li: q.li@xmu.edu.cn

    Yancheng You: yancheng.you@xmu.edu.cn

*Corresponding author：Qin Li; Email address: qin-li@vip.tom.com, q.li@xmu.edu.cn; Mailing address: School of Aerospace Engineering, Xiamen University, Xiang'an district, Xiamen, 361102, China

**Declarations**

- Availability of data and materials

  All data generated or analyzed during this study are included in this published article.

- Competing interests

  The authors declare that they have no competing interests.

- Funding

  Project of the National Numerical Wind-tunnel of China under grant number NNW2019ZT4-B12.

- Authors' contributions

  XF, YP and LQ analyzed and provided data regarding the quasi-linear analysis of dispersion relation preservation for nonlinear schemes; YY made substantial revision on the manuscript. All authors read and approved the final manuscript.

- Acknowledgements

  Thanks for Jianqiang Chen for his helpful discussion on the paper.

# Quasi-linear analysis of dispersion relation preservation for nonlinear schemes

Fengyuan Xu, Pan Yan, Qin Li, Yancheng You

School of Aerospace Engineering, Xiamen University, Xiamen, 361102, China

**Abstract:** In numerical simulations of complex flows with discontinuities, it is necessary to use nonlinear schemes. The spectrum of the scheme used have a significant impact on the resolution and stability of the computation. Based on the approximate dispersion relation method, we combine the corresponding spectral property with the dispersion relation preservation proposed by De and Eswaran (J. Comput. Phys. 218 (2006) 398–416) and propose a quasi-linear dispersion relation preservation (QL-DRP) analysis method, through which the group velocity of the nonlinear scheme can be determined. In particular, we derive the group velocity property when a high-order Runge–Kutta scheme is used and compare the performance of different time schemes with QL-DRP. The rationality of the QL-DRP method is verified by a numerical simulation and the discrete Fourier transform method. To further evaluate the performance of a nonlinear scheme in finding the group velocity, new hyperbolic equations are designed. The validity of QL-DRP and the group velocity preservation of several schemes are investigated using two examples of the equation for one-dimensional wave propagation and the new hyperbolic equations. The results show that the QL-DRP method integrated with high-order time schemes can determine the group velocity for nonlinear schemes and evaluate their performance reasonably and efficiently.

## 1 Introduction

Nonlinear approaches, such as the weighted essentially non-oscillatory (WENO) scheme, are widely used to calculate complex flow fields with discontinuities, such as shock waves. Because of the mathematics of such schemes, truncation errors can affect the convergence rate of the computation and the spectral properties affect the deviation of the Fourier modes of the numerical results from the exact solution. Therefore, studying the spectral properties of nonlinear schemes is important and meaningful when developing high-order nonlinear schemes.

Tam and Webb [1] suggested that it was necessary to consider the spectral properties besides using the standard Taylor series method in the construction of difference schemes. Specifically, by using the Fourier transformation of the governing equation, they proposed a basic method for analyzing the dispersion and dissipation relation of a linear scheme and developed a dispersion relation preservation (DRP) scheme. Moreover, this method can also be used to evaluate the spectral properties of different linear schemes and to construct new schemes. However, the method is inappropriate for analyzing the spectral properties of nonlinear schemes. To remedy this deficiency, Pirozzoli [2] studied the amplitude evolution of a disturbing mode after a tiny propagation period through a discrete Fourier transform (DFT), through which he proposed an approximate dispersion relation (ADR) analysis. ADR [2] gives an estimate of the total error generated by a nonlinear scheme and provides a method for analyzing its spectral properties. However, there are still some issues. For example, the time scheme can introduce errors, and if the grid number is not chosen carefully, there can be unexpected jumps in the spectral distribution [3].

Observing these problems, Mao et al. [3] proposed an improved ADR without time discretization called ADR-NT. Using the assumption of a tiny time period, they transferred the effect of the temporal derivative into the contribution of the spatial derivative. The researchers reported [3] that ADR-NT not only reduced the computational cost but also avoided the error due to the time discretization, as occurs in ADR. Mao et al. [3] also suggested ways to avoid jumping points when using ADR, i.e., the size of the grid should be twice some large prime or take the average over a large number of results with different initial phases.

The dispersion relation usually reflects the phase velocity of a spatial scheme. For practical problems with wave propagation, it is necessary to consider the group velocity, which is the one the energy propagates at [5]. De and Eswaran [4] analyzed DRP from the perspective of the group velocity, which we refer to as the DRP-g method or just DRP-g. They pointed out that DRP-g was important when evaluating the spectral properties of linear schemes. In the space of the reduced wave number and frequency, group velocity preservation (GVP) occurs in the region where the ratio of the numerical group velocity to the theoretical group velocity is in [0.95, 1.05]. Obviously, the analyses of DRP-g proposed by De and Eswaran [4] are different from those of Tam and Webb [1] and Pirozzoli [2]. However, the DRP-g method can be used only for linear schemes and not for nonlinear schemes. Thus, we derived the modified wave number of nonlinear schemes with the ADR method, combined the aforementioned group velocity analysis with high-order Runge–Kutta schemes, and finally, derived a quasi-linear method to analyze the spectral properties of nonlinear schemes. For brevity, the method is called the QL-DRP method or just QL-DRP. Its rationality was then verified with a numerical test and the DFT method.

As usual, a one-dimensional wave propagation equation [4] is used in this study to explore the GVP of nonlinear schemes by QL-DRP. As shown in [4], for this governing equation, the group velocity is the same as the phase velocity. Therefore, it is insufficient to investigate the GVP of difference schemes. To overcome this, we devised hyperbolic equations such that the solutions are synthetic waves where the group velocity and phase velocity are different. Using the equations, we numerically analyzed the GVP of different difference schemes and also verified the validity of QL-DRP.

This paper is arranged as follows. The DRP and ADR methods are reviewed in Section 2. In Section 3, the DRP-g formula is derived for a high-order Runge–Kutta scheme, the QL-DRP method is described, and the comparative results from QL-DRP are given for typical schemes. In Section 4, the group velocities of selected cases are obtained computationally using DFT and compared with those from QL-DRP, through which the rationality of QL-DRP is verified. We construct hyperbolic equations for which the solutions are waves where the group velocity and phase velocity are different. In Section 5, the one-dimensional wave propagation equation and the aforementioned hyperbolic equations are solved and analyzed. The conclusions are drawn in Section 6.

## 2 Review of DRP-g and ADR methods

Before we consider the QL-DRP method, the DRP-g and ADR methods are reviewed first.

### 2.1 DRP-g method

Consider the one-dimensional wave propagation equation [1–4]:

$$u_t + cu_x = 0. \tag{1}$$

For the initial distribution $u(x,0) = \hat{u}e^{ikx}$, Eq. (1) has the exact solution $u(x,t) = \hat{u}e^{i(kx-\omega t)}$. For Eq. (1) it holds that $\omega = kc$. Therefore, the exact group velocity is $V_{g,exact} = d\omega/dk = c$.

When using a time scheme, the following approximate relation is obtained:

$$u(t_n + \Delta t) - u^n \approx u^{n+1} - u^n, \tag{2}$$

where $u^n$ is the value at the initial time $t_n$, $u(t_n + \Delta t)$ is the exact solution at $t_n + \Delta t$, and $u^{n+1}$ is the corresponding numerical solution. Suppose the variable distribution on $t_n$ takes the form $u^n = \hat{u}e^{i(kx-\omega t)}$. When the explicit Euler scheme is used, Eq. (2) can be written as:

$$\hat{u}e^{i[kx-\omega(t+\Delta t)]} - \hat{u}e^{i(kx-\omega t)} \approx \Delta t \cdot f(t_n, u^n), \tag{3}$$

where $f(t_n, u^n)$ denotes the contributions from the spatial derivative. When a spatial scheme is used, one can see [4] that

$$f(t_n, u^n) = (-c)u_x = (-c)\frac{i\kappa'}{\Delta x}\hat{u}e^{i(kx-\omega t)},$$

where $\kappa'$ represents the modified wave number. It can be further derived from Eq. (3) that:

$$e^{-i\omega\Delta t} - 1 \approx \Delta t \cdot (-c) \cdot \frac{i\kappa'}{\Delta x}.$$

Following the approach of De and Eswaran [4], the numerical value $V_{g,num}$ of the Euler scheme can be derived by differentiating both sides of the above equation:

$$V_{g,num} = \left(\frac{d\omega}{dk}\right)_{num} = c\,\mathrm{Re}\left(e^{i\omega\Delta t}\cdot\frac{d\kappa'}{d\kappa}\right) = c\left(\cos(\omega\Delta t)\frac{d\kappa'_r(\kappa)}{d\kappa} - \sin(\omega\Delta t)\frac{d\kappa'_i(\kappa)}{d\kappa}\right), \quad (4)$$

where $\omega\Delta t$ is the reduced frequency, $\kappa = k\Delta x$ is the reduced wave number, and $c$ is the theoretical group/phase velocity. Equation (4) is exactly the same as that proposed by De and Eswaran [4].

When a linear scheme is employed, the spatial derivative at $x_i$ can be approximated by:

$$u_x(x_i) = \frac{1}{\Delta x}\sum_{j=-N}^{M} a_j u_{i+j},$$

where $M$ and $N$ are the numbers of nodes on the right and left sides of $x_i$, and $a_j$ are the scheme coefficients. The modified wave number can be obtained by a Fourier transformation [1]:

$$\kappa'(\kappa) = -i\sum_{j=-N}^{M} a_j e^{ij\kappa}. \quad (5)$$

Once the dispersion and dissipation relations are obtained, $V_{g,num}$ for the first-order explicit Euler scheme can be obtained from Eq. (4). The explicit Euler scheme yields relatively large temporal errors, and therefore, the group velocities obtained by DRP-g are inappropriate for unsteady problems when high-order Runge–Kutta schemes are used. Thus, Sengupta et al. [7] assumed that $a_4 + ib_4 = i\kappa'$ and $a_5 + ib_5 = \frac{d\kappa'}{d\kappa}$, and derived the following DRP-g formula for the fourth-order Runge–Kutta scheme:

$$\frac{V_{g,num}}{c} = \frac{a_5 + 2\sigma(b_4 b_5 - a_4 a_5)/3 + a_5\sigma^2(3a_4^2 - b_4^2 - N_c a_4^3 + N_c a_4 b_4^2)/6 + b_4\sigma^2(-6a_4 b_5 - 2b_4 a_5 - 3\sigma a_4^2 b_5 - \sigma b_5 b_4^2 + 2\sigma a_4 a_5 b_4)/6}{\cos(\omega\Delta t)}$$

where $\sigma = c\Delta t/\Delta x$ is the Courant–Friedrichs–Lewy (CFL) number. On checking, the following problems were found with this work:

1. There are no details for the derivation. Therefore, it is difficult to check the correctness of the formula.
2. When $\sigma$ approaches 0, the formula does not regress to that of the first-order explicit Euler scheme in Eq. (4).
3. When $\omega\Delta t = \pi/2$, the denominator $\cos(\omega\Delta t)$ is 0, and there is a singularity. Therefore, the DRP-g formula needs further investigation for higher-order Runge–Kutta schemes.

Thus, we employ the idea of [4] and derive the DRP-g results for high-order Runge–Kutta schemes in

Section 3.1.

## 2.2 ADR method

Despite its applicability to linear schemes, the spectral analysis method by Tam and Webb [1] is unsuitable for nonlinear schemes. Thus, Pirozzoli et al. [2] proposed a quasi-linear spectral analysis method for nonlinear schemes or ADR. The process is as follows [2].

For Eq. (1), suppose the initial distribution is $u(x,0) = \hat{u}_0 e^{ikx}$ and consider a semi-discretized approximation of Eq. (1) on grids $\{x_j = j\Delta x\}$ with the spacing $\Delta x$ as

$$\frac{dv_j}{dt} + c\delta v_j = 0, v_j(0) = \hat{u}_0 e^{ij\kappa}, \qquad (6)$$

where $v_j(t) \approx u(x_j, t)$ and $\delta v_j$ is for a specific numerical scheme, e.g., the linear difference scheme:

$$\delta v_j = \frac{1}{\Delta x} \sum_{l=-N}^{M} a_l u_{j+l},$$

as before. From [2], Eq. (6) has the following exact solution: $v_j(t) = \hat{v}(t) e^{ij\kappa} = \hat{u}_0 e^{-i(ct/\Delta x)\kappa'} e^{ij\kappa}$. Furthermore, Pirozzoli [2] proposed ADR to derive $\kappa'$:

$$\kappa'(\kappa_n) = \frac{i\Delta x}{c\tau} \ln\left(\frac{\hat{v}(\kappa_n, \tau)}{\hat{v}_0(\kappa_n)}\right). \qquad (7)$$

In this equation, $\hat{v}_0(\kappa_n)$ is the Fourier coefficient corresponding to the initial mode and $\hat{v}(\kappa_n, \tau)$ is the DFT of $v_j(\tau)$ at time $\tau$ where $\kappa_n = k_n \Delta x$. The operation of DFT is

$$\hat{v}(\kappa_n, \tau) = \frac{1}{N_x} \sum_{j=0}^{N_x - 1} v_j(\tau) e^{-ij\kappa_n}, \qquad (8)$$

where $N_x$ is the number of grid number.

As indicated in the introduction, ADR depends on the choice of time step and number, and inappropriate choices can lead to temporal errors or unreasonable jump points in the spectral distributions [3]. Mao et al. [3] expanded the solution $v_j(\tau)$ of Eq. (6) with a Taylor series:

$$v_j(\tau) = v_j(0) + \tau\left[-\frac{c}{\Delta x}\sum_{l=-N}^{M} b_{j+l} v_{j+l}(0)\right] + o(\tau^2).$$

Substituting the above formula into Eq. (8) and the subsequent result into Eq. (7), Mao et al. [3] proposed ADR without time discretization, known as ADR-NT, to derive $\kappa'(\kappa_n)$:

$$\kappa'(\kappa_n) = \frac{1}{i}\left(\sum_{j=0}^{N_x-1}\left(\sum_{l=-N}^{M} b_{j+l} v_{j+l}(0)\right) e^{-ij\kappa_n}\right) \Big/ (N_x \widehat{v}_0), \tag{9}$$

where $b_{j+l}$ are the coefficients of the nonlinear scheme. Compared with ADR, ADR-NT reduces the computation cost and minimizes the influence of temporal errors.

## 3  QL-DRP for nonlinear schemes

De and Eswaran [4] first analyzed and derived the numerical group velocity for linear schemes. However, their method can be used only for linear schemes and not for nonlinear ones. Considering the development of ADR [2] from DRP [1], which extends the analysis from linear to nonlinear, we introduce the influence of the nonlinearity in DRP-g [4] by using ADR to account for the spectral property of the nonlinear scheme. Next, employing a similar procedure as in DRP-g, the so-called QL-DRP is derived for nonlinear schemes. Prior to a further discussion on QL-DRP, we first derive the DRP-g formula for high-order Runge–Kutta schemes.

### 3.1  DRP-g formula for high-order Runge–Kutta schemes

Taking the third-order Runge–Kutta scheme as an example, we derive the DRP-g formula based on the method of [4]. For convenience, the $n$th-order Runge–Kutta scheme is referred to as $RK_n$.

$RK_3$ here takes the form [8]:

$$\begin{cases} u^{(1)} = u^n + \Delta t \cdot f(u^n), \\ u^{(2)} = \dfrac{3}{4} u^n + \dfrac{1}{4} u^{(1)} + \dfrac{1}{4}\Delta t \cdot f(u^{(1)}), \\ u^{n+1} = \dfrac{1}{3} u^n + \dfrac{2}{3} u^{(2)} + \dfrac{2}{3}\Delta t \cdot f(u^{(2)}), \end{cases} \tag{10}$$

where $f(u^n)$, $f(u^{(1)})$, and $f(u^{(2)})$ are the derivatives on the sub-time steps. From Eq. (10), $u^{n+1}$ can be reformulated as:

$$u^{n+1} = u^n + \Delta t\left(\frac{1}{6} f(u^{(n)}) + \frac{1}{6} f(u^{(1)}) + \frac{2}{3} f(u^{(2)})\right).$$

Let $K_1 = f(u^{(n)})$, $K_2 = f(u^{(1)})$, and $K_3 = f(u^{(2)})$. Then, Eq. (2) can be rewritten as:

$$\widehat{u} e^{i[kx-\omega(t+\Delta t)]} - \widehat{u} e^{i(kx-\omega t)} \approx \Delta t \cdot \left(\frac{1}{6} K_1 + \frac{1}{6} K_2 + \frac{2}{3} K_3\right). \tag{11}$$

It is known that $K_1$ in Eq. (10) corresponds to the derivative of $u^n$ or $K_1 = f(t^n, u^n)$; similarly

$K_2 = f(t^n + \Delta t, u^{(1)}) = f\left(t^n + \Delta t, u^n + \Delta t K_1\right)$ and $K_3 = f\left(t^n + \tfrac{1}{2}\Delta t, u^n + \tfrac{1}{2}\Delta t \cdot \tfrac{1}{2}(K_1 + K_2)\right)$. Using $\kappa'$, then $K_i$ in Eq. (11) can be written as [4,9]:

$$\begin{cases} K_1 = (-c) \cdot \dfrac{i\kappa'}{\Delta x} \hat{u} e^{i(kx-\omega t)}, \\ K_2 = \left[(-c) \cdot \dfrac{i\kappa'}{\Delta x} + \Delta t \cdot (-c^2)\left(\dfrac{\kappa'}{\Delta x}\right)^2\right] \hat{u} e^{i(kx-\omega t)}, \\ K_3 = \left[(-c) \cdot \dfrac{i\kappa'}{\Delta x} + \dfrac{\Delta t}{4} \cdot (-c^2)\left(\dfrac{\kappa'}{\Delta x}\right)^2 + \dfrac{\Delta t}{4}\left((-c^2)\left(\dfrac{\kappa'}{\Delta x}\right)^2 + \Delta t \cdot ic^3\left(\dfrac{\kappa'}{\Delta x}\right)^3\right)\right] \hat{u} e^{i(kx-\omega t)}. \end{cases} \quad (12)$$

Substituting Eq. (12) into Eq. (11):

$$e^{-i\omega\Delta t} - 1 \approx -i\sigma \cdot \kappa' - \dfrac{1}{2}(\sigma \cdot \kappa')^2 + \dfrac{1}{6}i(\sigma \cdot \kappa')^3.$$

Applying $d(.)/dk$ to the above equation and taking the real part as in [4], the group velocity for $RK_3$ is

$$\left(V_{g,num}\right)_{RK_3} = c \cdot \mathrm{Re}\left((1 - i\sigma \cdot \kappa' - \dfrac{1}{2}(\sigma \cdot \kappa')^2) e^{i\omega\Delta t} \cdot \dfrac{d\kappa'}{d\kappa}\right). \quad (13)$$

For a linear scheme, $\kappa'$ is found from Eq. (5) [4]. Similarly, when an $RK_4$ scheme is used for time discretization, the DRP-g formula can be derived as:

$$\left(V_{g,num}\right)_{RK_4} = c \cdot \mathrm{Re}\left((1 - i\sigma \cdot \kappa' - \dfrac{1}{2}(\sigma \cdot \kappa')^2 + \dfrac{1}{6}i(\sigma \cdot \kappa')^3) e^{i\omega\Delta t} \cdot \dfrac{d\kappa'}{d\kappa}\right). \quad (14)$$

From Eqs. (13) and (14), it can be seen that the numerical group velocity for a high-order Runge–Kutta scheme depends on $\sigma$, unlike the first-order explicit Euler scheme in Eq. (4). For illustration, the group velocity of a linear fifth-order scheme, UPW5, namely:

$$(du/dx)_j \approx \tfrac{1}{\Delta x}\left(-\tfrac{1}{30}u_{i-3} + \tfrac{1}{4}u_{i-2} - u_{i-1} + \tfrac{1}{3}u_i + \tfrac{1}{2}u_{i+1} - \tfrac{1}{20}u_{i+2}\right),$$

was investigated for different CFL numbers. The velocity distributions are shown in Figs. 1 and 2 for $\sigma = 0.01$ and 0.1, respectively. The effect of $\sigma$ is apparent through the variation of the GVP region, for which $V_{g,num}/c \in [0.95, 1.05]$, as mentioned before.

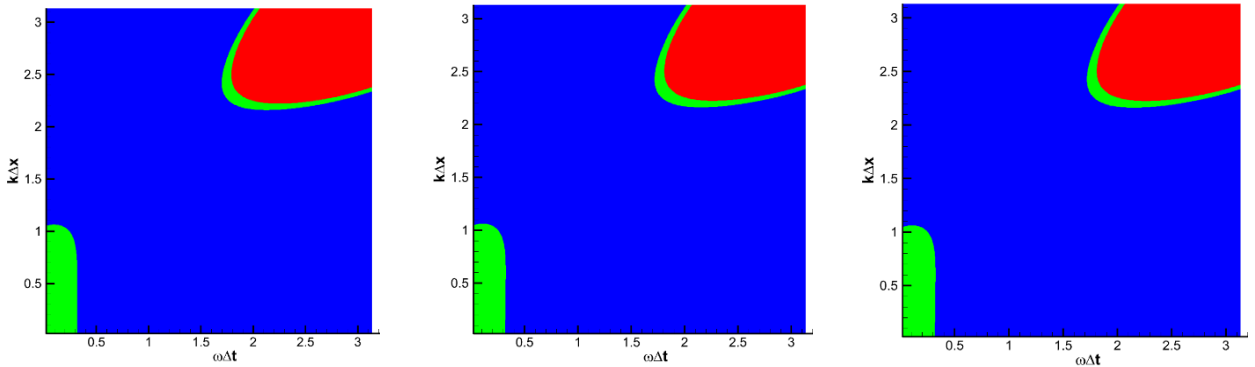

(a) Explicit Euler scheme　　　　　　　　(b) $RK_3$　　　　　　　　(c) $RK_4$

Fig. 1. Distribution of $v_g/c$ in the $k\Delta x$ vs. $\omega\Delta t$ space for different time schemes. UPW5 is used for the spatial discretization, and $\sigma = 0.01$. The green region represents the range 0.95–1.05, the red region represents the range >1.05, and the blue region represents the range <0.95.

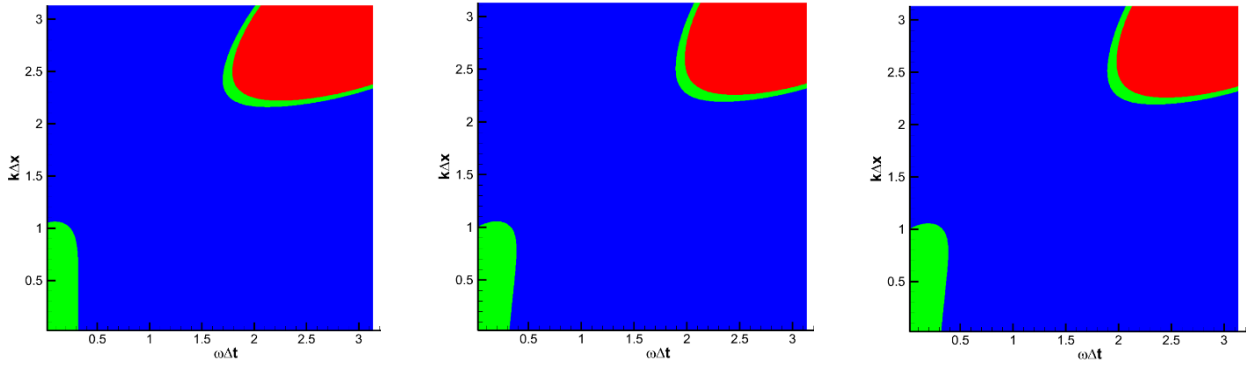

(a) Explicit Euler scheme　　　　　　　　(b) $RK_3$　　　　　　　　(c) $RK_4$

Fig. 2. Distribution of $v_g/c$ in the $k\Delta x$ vs. $\omega\Delta t$ space for different time schemes. UPW5 is used for the spatial discretization, and $\sigma = 0.1$. The green region represents the range 0.95–1.05, the red region represents the range >1.05, and the blue region represents the range <0.95.

As indicated by De and Eswaran [4], the group distribution near the origin is of special concern. Figure 1 shows that when the CFL number is small, the distributions of the first-order explicit Euler scheme, $RK_3$, and $RK_4$ are almost the same. With an increase of the CFL number, differences appear among the different time schemes. As shown in Fig. 2, when $\sigma = 0.1$, for the green area near the origin, $RK_4 \approx RK_3 >$ the explicit Euler scheme, which indicates that higher-order time schemes enhance the GVP. This outcome confirms the potential advantages of higher-order time schemes for unsteady problems.

In summary, in this section the group velocities of $RK_3$ and $RK_4$ were found by applying the same procedure used for the explicit Euler scheme in [4]. The following remarks are given: (1) When the CFL number approaches 0, the group velocity distributions of higher-order Runge–Kutta schemes gradually degenerate into that of the first-order explicit Euler scheme. Further, the results from DRP-g for different time schemes become the same under this condition. (2) No singularity arises in Eqs. (13) and (14), which indicates their rationality.

### 3.2　QL-DRP

As mentioned in the introduction, QL-DRP is proposed to evaluate the characteristics of the group velocity in a quasi-linear manner. The main difference between QL-DRP and DRP-g is that the modified wave number in the former considers the nonlinearity of schemes. Therefore, the formulas for QL-DRP are the same as those for DRP-g, such as Eqs. (13) or (14), except that different evaluations of $\kappa'$ are implemented.

The concrete implementations are summarized as follows:

1. Solve for the modified wave number $\kappa'(\kappa)$ of nonlinear spatial schemes.

Based on ADR, the dispersion and dissipation relations of the nonlinear scheme, $\kappa'(\kappa)$, are solved with Eq. (9).

2. Solve for $d\kappa'/d\kappa$.

When the spatial scheme is linear, $\kappa'(\kappa)$ can be obtained as $\kappa'(\kappa) = -i \sum_{j=-N}^{M} a_j e^{ij\kappa}$. Because the coefficient $a_j$ is constant, $d\kappa'/d\kappa$ can be derived analytically. However, when the spatial scheme is nonlinear, the coefficient $b_{j+l}$ in Eq. (9) correlates with the initial variable distribution nonlinearly and $d\kappa'/d\kappa$ is hard to resolve analytically. To overcome this difficulty, we use the difference method to evaluate it numerically. Taking the second-order central difference as an example, $d\kappa'/d\kappa$ can be evaluated as:

$$\left(\frac{d\kappa'}{d\kappa}\right)_{j+1/2} \approx \frac{\kappa'_{j+1} - \kappa'_j}{\Delta \kappa}, \tag{15}$$

where $j$ is the index of the reduced wave number and

$$\kappa_j = k_j \Delta x = \frac{2\pi}{\lambda_j} \times \frac{L}{j} = \frac{2\pi j}{N_x},$$

for $j = 0, \ldots, N_x/2$.

3. Substitute $d\kappa'/d\kappa$ into the QL-DRP formulas, such as Eqs. (13) or (14), to obtain $V_{g,num}$.

When solving for $\kappa'(\kappa)$ for nonlinear schemes with ADR, $\Delta t$ in the time integration must be small enough [2], which means the corresponding reduced frequency $\omega \Delta t$ should also be small. Correspondingly, QL-DRP can describe the GVP of nonlinear schemes effectively only for a small reduced frequency. When $\omega \Delta t$ is large, ADR does not accurately predict the spectral property, and therefore, QL-DRP cannot accurately describe GVP either. However, QL-DRP can give some insights into the group velocity of a nonlinear scheme. Hence, QL-DRP can provide a useful reference for the overall characteristics of the group velocity for a nonlinear scheme.

For demonstration, we used QL-DRP to find the group velocity for three typical spatial schemes: UPW5, WENO5-JS [8], and WENO5-M [6]. Adopting the aforementioned procedure, the corresponding distributions of $V_g/c$ with the $RK_4$ scheme are shown in Fig. 3.

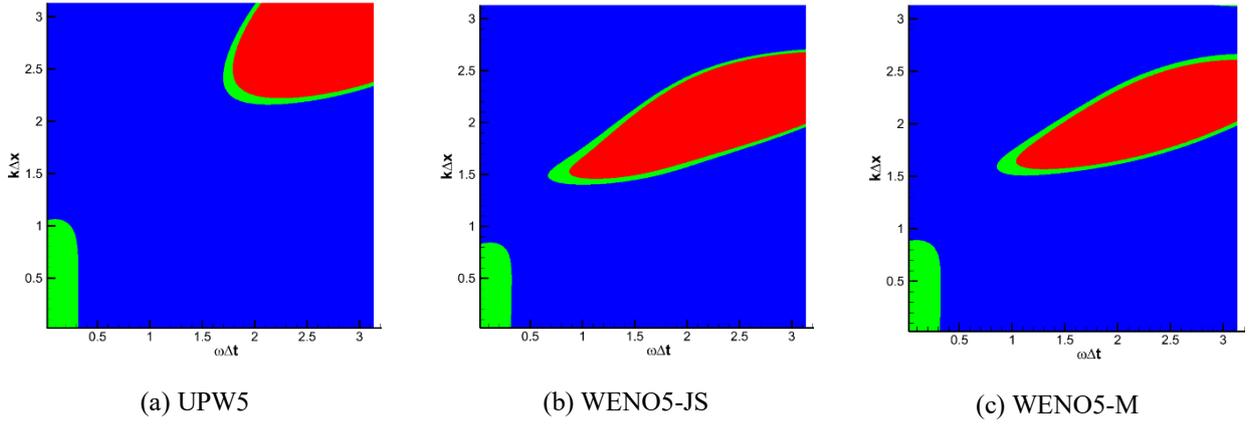

(a) UPW5  (b) WENO5-JS  (c) WENO5-M

Fig. 3. Distribution of $v_g/c$ in the $k\Delta x$ vs. $\omega\Delta t$ space for different spatial schemes. $RK_4$ is used for the time discretization, and $\sigma = 0.01$. The green region represents the range 0.95–1.05, the red region represents the range >1.05, and the blue region represents the range <0.95.

The figure shows that:

1. Comparing Fig. 3(a) with Fig. 1, one can see that the group velocity distributions for QL-DRP with the UPW5 scheme are the same as those for DRP-g. This is because, for linear schemes, $\kappa'(\kappa)$ in ADR is the same as that in DRP [1]. In other words, QL-DRP degenerates into DRP-g for a linear spatial scheme.

2. The distributions of the group velocity are obviously different for the three schemes, which indicates that QL-DRP can be used to distinguish between the GVP of different spatial schemes under the same time scheme. Specifically, the area of the GVP region (green) for the WENO5-JS scheme is smaller than that for the UPW5 scheme, which indicates that the group velocities of the nonlinear scheme differ from that of linear scheme. Moreover, the GVP area of the WENO5-M scheme is larger than that of WENO5-JS, which shows that the nonlinear optimization of the former enhances the GVP.

3. The GVP area near the origin, especially along the vertical axis, is the largest for UPW5, followed by WENO5-M and then WENO5-JS.

In short, QL-DRP can be applied not only to analyze the group velocity of a linear scheme but also of a nonlinear scheme. Moreover, the method can clearly differentiate spatial schemes by considering the size of the GVP region.

## 4  Numerical analysis and validation of QL-DRP

The previous section introduced the derivation and implementation of QL-DRP. In this section, the numerical group velocity is obtained numerically and the results are compared with those from QL-DRP. Thus, the rationality of the method was analyzed and verified. Moreover, we devise hyperbolic equations and a corresponding example to determine the numerical group velocity of the scheme directly, which further verifies the approach.

**4.1 Numerical analysis and validation**

In this section, the group velocity of the different schemes is numerically obtained by solving the one-dimensional wave propagation equation and by using DFT. By comparing with the analytic results from QL-DRP, the numerical error is acquired and the rationality of QL-DRP is evaluated.

First, we select several points $P_i$ in the $k\Delta x$ vs. $\omega\Delta t$ plane, then numerically compute the group velocity using DFT. The coordinates of $P_i$ are the reduced frequency and wave number. Without loss of generality, we set $k\Delta x$ as 1. Because $\omega = kc$, then

$$\omega\Delta t = c\frac{\Delta t}{\Delta x}\kappa.$$

Hence, once $\Delta t$ and $\Delta x$ are determined, the coordinates of $P_i$ are defined.

Next, the derivation of numerical group velocity will be discussed. The solution of the semi-discrete equation (6) is

$$v_j(t) = \hat{u}_0 e^{-i(ct/\Delta x)\kappa'} e^{ij\kappa} = \hat{u}_0 \exp\left[i\left(kx - \frac{c\kappa'}{\Delta x}t\right)\right]$$

where $\kappa' = \kappa'(\kappa)$ [10]. Compared with the exact solution $u_j(t) = \hat{u}_0 e^{i(kx-\omega t)}$, the numerical frequency can be expressed as $\omega' = c\kappa'/\Delta x$, which depends on $\kappa'$. Hence, the reduced frequency at $P_i$ can be expressed with the modified wave number as:

$$\omega'\Delta t = c\frac{\kappa'\Delta t}{\Delta x} = c\kappa'\frac{N\Delta t}{L}, \tag{16}$$

where $L$ is the length of the computational domain. We consider Eq. (1) with $c = 1$ and the initial distribution $u_0 = \cos(x)$ at $x \in [0, L]$ with $L = 2\pi$. From [2], $\kappa'$ can be numerically found with ADR [2] for a spatial scheme and time scheme using DFT. Here, we use WENO5-JS and $RK_4$ for the spatial and time discretization. The number of grid cells and time step are given in Table 1. As shown in the table, four cases with respective ($N_x$, $\Delta t$) are chosen, and these correspond to $P_i$. Then, $\kappa'$ can be acquired with respect to $\kappa$, as illustrated in Fig. 4.

Table 1. Size of computational grid, time step, and coordinates of $P_i$.

| Point | Grid number $N_x$ | Time step, $\Delta t$ | $(k\Delta x, \omega\Delta t)$ |
|---|---|---|---|
| $P_1$ | 422 | $10^{-8}$ | $(1, 6.72\times 10^{-7})$ |
| $P_2$ | 422 | $10^{-3}$ | (1, 0.0672) |
| $P_3$ | 3046 | $10^{-3}$ | (1, 0.04848) |
| $P_4$ | 6082 | $10^{-3}$ | (1, 0.968) |

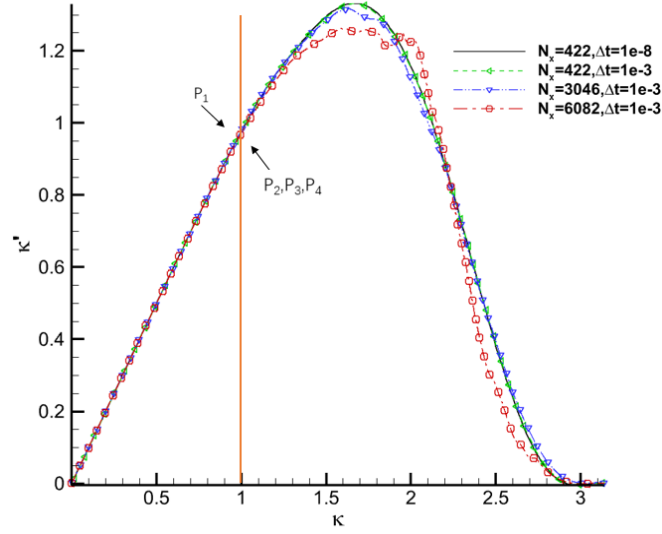

Fig. 4. Illustrations of $P_i$ on the distributions of $\text{Re}(\kappa')$ for different sizes of grid and different time steps with WENO5-JS and $RK_4$.

Besides $P_1$ to $P_4$, we also investigated the case when the reduced frequency was larger. There were oscillations in the results. Hence, the frequency region for investigation was limited to $\omega \Delta t < 1$. For illustration, the positions of $P_i$ ($i$ = 1, 2, 3, 4) are shown in the $k\Delta x$ vs. $\omega \Delta t$ plane in Fig. 5. The colors indicate the group velocity calculated by QL-DRP with WENO5-JS and $RK_4$.

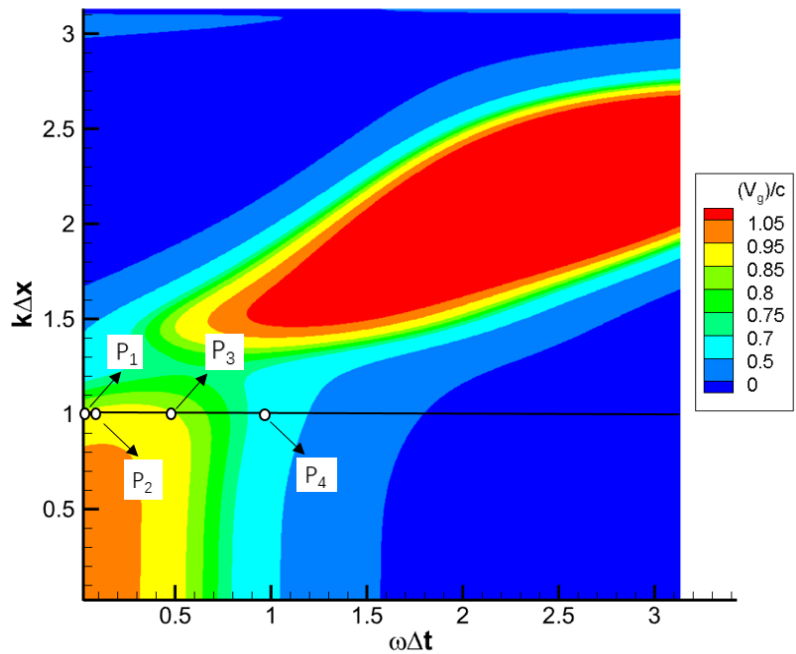

Fig. 5. Positions of $P_i$ in the $k\Delta x$ vs. $\omega \Delta t$ plane. The colors indicate the group velocity calculated by QL-DRP with

Thus far, the numerical solution of the group velocity can be found from:

$$V_{g,num} = \mathrm{Re}(\Delta\omega')/\Delta k = \mathrm{Re}(\omega'_2 - \omega'_1)/(k_2 - k_1), \tag{17}$$

where $k_1$ and $k_2$ represent two wave numbers with a small difference, $\omega'_1$ and $\omega'_2$ are the corresponding numerical frequencies, and Re denotes taking the real part. Substituting Eq. (16) into Eq. (17):

$$V_{g,num} = \frac{\mathrm{Re}(\kappa'_2(\kappa_2) - \kappa'_1(\kappa_1))}{\kappa_2 - \kappa_1}. \tag{18}$$

Under the assumption that $\kappa = 1$ as above, we set $\kappa_1$ and $\kappa_2$ as $\kappa_1 = 1 - \Delta\kappa$ and $\kappa_2 = 1 + \Delta\kappa$. The corresponding modified wave numbers $\mathrm{Re}(\kappa'_1(\kappa_1))$ and $\mathrm{Re}(\kappa'_2(\kappa_2))$ can be obtained from the distributions of $\mathrm{Re}(\kappa')$ in Fig. 4. Moreover, $V_{g,num}$ is the slope of $\mathrm{Re}(\kappa')$ at the point $(1, \kappa'_i(1))$.

Next, we compare the group velocities at $P_i$ obtained by QL-DRP and the numerical analysis. For QL-DRP, the analytic group velocity $V_{g,anal}$ at $P_i$ can be obtained from Eqs. (9) and (14), whereas the numerical group velocity $V_{g,num}$ at $P_i$ can be obtained from Eq. (18). These are compared in Table 2.

Table 2. Both group velocities and errors of $P_i$ for WENO5-JS and $RK_4$.

| Point | $V_{g,num}$ | $V_{g,anal}$ | $\lvert V_{g,anal} - V_{g,num}\rvert / V_{g,num}$ |
|---|---|---|---|
| $P_1$ | 0.8647 | 0.8627 | 0.23% |
| $P_2$ | 0.8638 | 0.8698 | 0.69% |
| $P_3$ | 0.8732 | 0.8524 | 2.38% |
| $P_4$ | 0.8203 | 0.6505 | 20.70% |

Table 2 shows that when $\omega\Delta t$ is small, as for $P_1$, $P_2$, and $P_3$, then the group velocity given by QL-DRP is quite close to the numerical value. This indicates that QL-DRP can reasonably predict the group velocity at low reduced frequency. With an increase of $\omega\Delta t$, the error between $V_{g,anal}$ and $V_{g,num}$ becomes large but was still within 21% when $\omega\Delta t < 1$. Therefore, QL-DRP can provide important insights for the group velocity at medium reduced frequency.

Note that QL-DRP can provide an overview of the group velocity in the $k\Delta x$ vs. $\omega\Delta t$ space, which is difficult for the numerical approach just mentioned. In addition, QL-DRP avoids the tedious process of finding a numerical solution and is more convenient for practical analyses.

## 4.2 Constructing hyperbolic equations for a combination wave with different phase and group velocities

The distribution $u = \cos(kx - \omega t)$ satisfies Eq. (1) when $\omega/k = c$, where the group and phase velocities are $V_{g,exact} = d\omega/dk = c$ and $V_{p,exact} = \omega/k = c$, respectively. One can see that in such situations, the group velocity and phase velocity are indistinguishable, which is unfavorable for investigating the GVP property of the numerical scheme.

In textbooks, the combination wave $u = \cos(k_1 x - \omega_1 t) + \cos(k_2 x - \omega_2 t)$ is usually used to explain the concept of group velocity, which can also be written as:

$$u = 2\cos\left(\frac{k_1 + k_2}{2} x - \frac{\omega_1 + \omega_2}{2} t\right) \cos\left(\frac{k_2 - k_1}{2} x - \frac{\omega_2 - \omega_1}{2} t\right). \tag{19}$$

The envelope of Eq. (19) is

$$2\cos\left(\frac{k_2 - k_1}{2} x - \frac{\omega_2 - \omega_1}{2} t\right),$$

and the corresponding group velocity and phase velocity are

$$V_{g,exact} = \frac{\omega_2 - \omega_1}{k_2 - k_1} \quad \text{and} \quad V_{p,exact} = \frac{\omega_2 + \omega_1}{k_2 + k_1}.$$

Without loss of generality, it is usually assumed that $k_1 = \omega_1$. When $k_2 \neq \omega_2$, then $V_g \neq V_p$.

Although the combination wave is used to explain the group velocity and phase velocity, the distribution does not satisfy Eq. (1). Therefore, this equation cannot be used to simulate the wave. We devised the following hyperbolic equations so that we could numerically study the GVP of different schemes:

$$\begin{cases} \dfrac{\partial u}{\partial t} + \dfrac{\partial u}{\partial x} = p, \\ \dfrac{\partial p}{\partial t} + a\dfrac{\partial p}{\partial x} = 0. \end{cases} \tag{20}$$

It can be shown that Eq. (19) satisfies Eq. (20) and $a = \omega_2/k_2$ when $k_1 = \omega_1$. If $a \neq 1$, the group velocity is different from the phase velocity. Hence, a measure is provided to study the combination wave with the co-occurrence of different group and phase velocities, which also provides a powerful way to verify the outcome of QL-DRP.

## 5 Numerical examples

In the following, to compute the group velocity, the envelope of the wave is derived by a Hilbert transform.

### 5.1 One-dimensional wave propagation

Consider Eq. (1) with $c = 1/8$. The following initial distribution is chosen:

$$u_0 = \sin(8\pi x). \tag{21}$$

The exact solution is $u = \sin(8\pi x - \pi t)$ with $V_{g,exact} = V_{p,exact} = 1/8$.

Two tests were used to check the GVP: (1) using the same spatial scheme for different $\kappa$ and (2) using different spatial schemes for the same $\kappa$. Considering that $\kappa = k\Delta x = kL/N_x$, $\kappa$ can be adjusted by changing the size of $N_x$. In the computation, $T = 2$, $\Delta t = 10^{-3}$, and $\omega\Delta t = 10^{-3}\pi$, and the computational domain was $[-1, 1]$.

First, the performance of WENO5-JS was tested for different grids, as shown in Fig. 6.

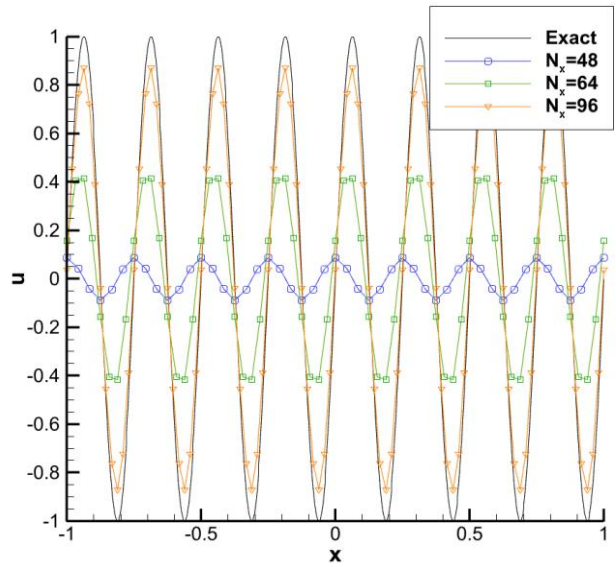

Fig. 6. Results of WENO5-JS with $RK_4$ for different sizes of grid for $T = 2$, $\Delta t = 10^{-3}$, and $\omega\Delta t = 10^{-3}\pi$.

When $N_x$ was 48, 64, or 96, $\kappa$ was $\pi/3$, $\pi/4$, or $\pi/6$, respectively. Because the group velocity is equal to the phase velocity in this situation, the phase velocity can be computed by the phase change of the wave peak, and the group velocity can be obtained thereafter. Figure 6 shows that when $\kappa = \pi/3$, the numerical group velocity was about 0.1025 and the ratio $V_{g,num}/V_{g,exact} = 0.82$. When $\kappa = \pi/4$, $V_{g,num}$ was about 0.11875 and $V_{g,num}/V_{g,exact} = 0.95$, which indicated that the group velocity was in the accurate region [0.95, 1.05]. With an increase of $N_x$, $\kappa$ decreased and the numerical group velocity gradually approached the theoretical value, i.e., $V_{g,num}$ was 0.99 when $\kappa = \pi/6$.

Using QL-DRP with the WENO5-JS scheme and the $RK_4$ scheme, for the same three points $(\pi/3, 10^{-3}\pi)$, $(\pi/4, 10^{-3}\pi)$, and $(\pi/6, 10^{-3}\pi)$ in the $k\Delta x$ vs. $\omega\Delta t$ plane, $V_{g,num}/V_{g,exact}$ can be derived from Eq. (14) as 0.8259, 0.9592, and 0.9950, respectively. These values are nearly the same as those obtained numerically. Hence, QL-DRP has been verified quantitatively, and its usefulness in analyses of the group velocity has been demonstrated.

Next, GVP was compared for UPW5, WENO5-JS, and WENO5-M for $N_x$ = 48 and $(\kappa, \omega\Delta t) = (\pi/3, 10^{-3}\pi)$. The results are shown in Fig. 7.

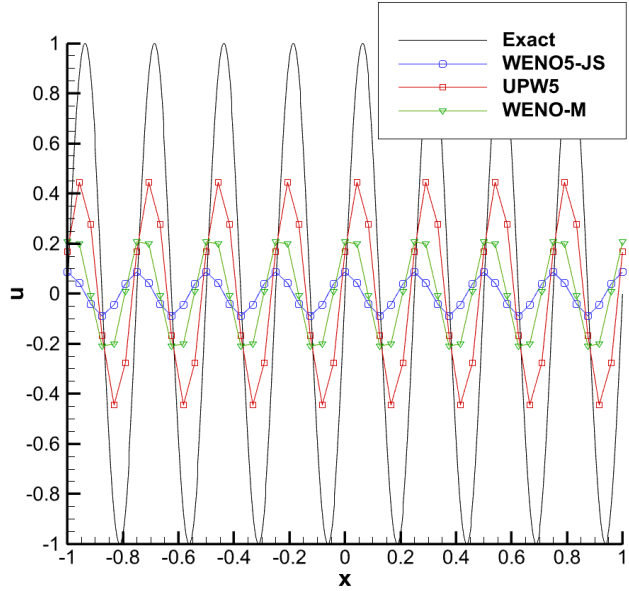

Fig. 7. Results for UPW5, WENO5-JS, and WENO5-M with $RK_4$ for grid size $N_x$ = 48, $T$ = 2, $\Delta t = 10^{-3}$, and $\omega\Delta t = 10^{-3}\pi$.

Figure 7 shows that the UPW5 scheme yields a relatively accurate group velocity, since $V_{g,num}/V_{g,exact} = 0.96$, whereas WENO5-M and WENO5-JS have relatively large errors, i.e., $V_{g,num}/V_{g,exact} = 0.87$ and 0.82, respectively. WENO5-M performed better than WENO5-JS. The rank of performance for the GVP is UPW5 > WENO5-M > WENO5-JS. The results for QL-DRP are consistent with the numerical results, which demonstrates the validity and capability of QL-DRP for nonlinear schemes.

## 5.2 Wave propagation with different group and phase velocities

In this section, the case described in Section 4.2 with different group and phase velocities is evaluated for the different spatial schemes, namely UPW5, WENO5-M, and WENO5-JS, which further demonstrates the validity of QL-DRP.

For the distribution in Eq. (19), f multiple waves were included in one wave packet in purpose of better illustration. The initial condition is

$$u_0 = 2\cos\left(\frac{k_1+k_2}{2}x\right)\cos\left(\frac{k_2-k_1}{2}x\right)$$

in the domain $[-3\pi, 3\pi]$. We set $k_1 = \omega_1 = 6$, $k_2 = 8$, and $\omega_2 = 12$. The exact group and phase velocities were $V_{g,exact} = 3$ and $V_{p,exact} = 9/7$. The initial distribution and its envelope derived with a Hilbert transform are shown in Fig. 8. Here, $N_x = 120$, $T = 1$, and $\Delta t = 5 \times 10^{-4}$.

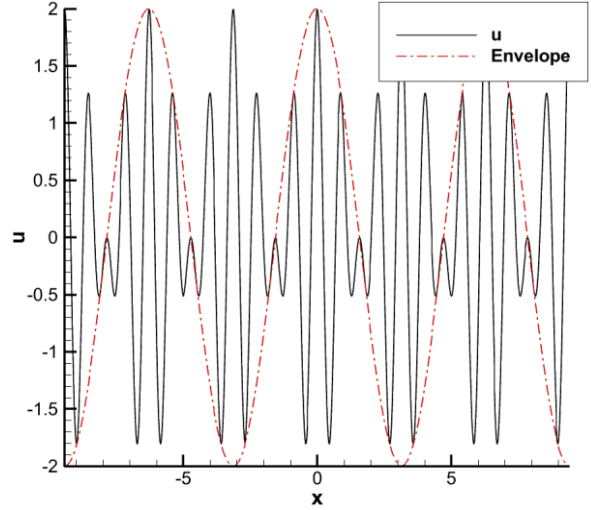

Fig. 8. Distributions of initial combination wave and its envelope for $k_1 = \omega_1 = 6$, $k_2 = 8$, and $\omega_2 = 12$.

The corresponding envelopes, which contain information about the group velocity, were derived, as shown in Fig. 9. By convention, the envelopes were drawn after taking the absolute value.

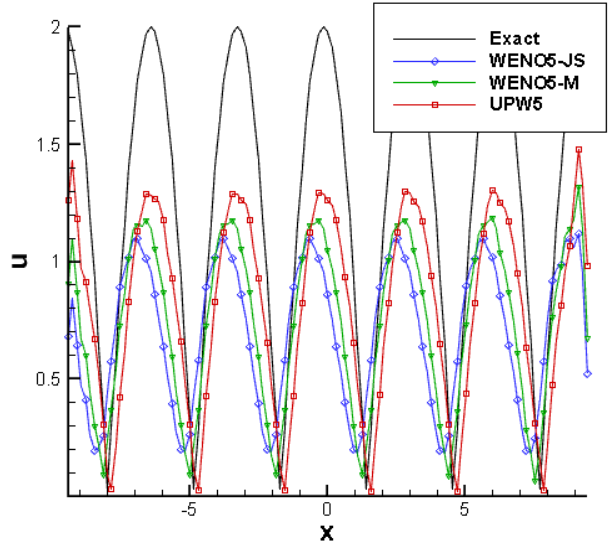

Fig. 9. Envelopes for $T = 1$ for the hyperbolic equations under the conditions $N_x = 120$ and $\Delta t = 5 \times 10^{-4}$ with $k_1 = \omega_1 = 6$, $k_2 = 8$, and $\omega_2 = 12$ for different spatial schemes and $RK_4$.

Figure 9 shows that UPW5 yields an envelope such that the numerical group velocity is closest to the theoretical one, followed by WENO5-M and then WENO5-JS. The amplitudes of the envelope show that UPW5 had the least dissipation, followed by WENO5-M and then WENO5-JS. Recall that from Fig. 3, the area of GVP given by QL-DRP was in the order: UPW5 > WENO5-M > WENO5-JS. Therefore, the results for QL-DRP agree with the numerical results, which verifies its validity.

# 6　Conclusions

Because DRP-g [4] cannot be used for a group velocity analysis of a nonlinear scheme, QL-DRP was proposed by combining ADR with DRP-g. Moreover, a detailed derivation and implementation of QL-DRP were provided for high-order Runge–Kutta schemes. The conclusions are as follows:

1. Since QL-DRP can analyze the GVP for nonlinear schemes, it can distinguish between the group velocity of different schemes. When the spatial scheme is linear, QL-DRP degenerates to DRP-g.
2. The group velocities of typical schemes were derived numerically and the comparison confirmed the rationality of QL-DRP.
3. If the reduced frequency is small, the group velocity given by QL-DRP is accurate. Moreover, its predictions are reasonable and meaningful at medium reduced frequency. Although QL-DRP cannot accurately predict the GVP at high reduced frequency, it can provide a reference in an analysis of the group velocity of nonlinear schemes.
4. Hyperbolic equations were devised for combination waves with different group and phase velocities. These enabled numerical investigations. Through two numerical examples, the validity of QL-DRP was verified qualitatively and quantitatively.

**Acknowledgements**

This study is sponsored by a project of the National Numerical Wind-tunnel of China under grant number NNW2019ZT4-B12.